\documentclass[11pt,a4paper,reqno]{amsart}
\usepackage{amsthm, amssymb, amscd} 
\usepackage{enumerate,array}
\usepackage{oppi}
\newtheorem{theorem}{Theorem}[section]
\newtheorem{lemma}[theorem]{Lemma}

\newtheorem{proposition}[theorem]{Proposition}

\newtheorem*{capu}{The $\a$-shuffles conjecture}
\theoremstyle{definition}

\theoremstyle{remark}
\newtheorem{remark}[theorem]{Remark}

\newtheorem*{HNA}{(Harmless) notation ambiguity}
\numberwithin{equation}{section}
\newcommand{\hw}{\widehat w}
\newcommand{\hX}{\widehat X}
\newcommand{\hY}{\widehat Y}

\begin{document} 
%------------------------------------------------------------------------
%

\title[Octopus inequality]{%
A few remarks on the octopus inequality and Aldous' spectral gap conjecture%
}
%    author one information
\author[Filippo Cesi]{Filippo Cesi}
\address{%
    Filippo Cesi\hfill\break
    \indent Dipartimento di Fisica\hfill\break
    \indent Universit\`a di Roma ``La Sapienza", Italy.
    %\hfil\break
    %\indent and SMC, INFM-CNR.
}
%\curraddr{}
\email{filippo.cesi@roma1.infn.it}
%\thanks{}

%    \subjclass is required by all journals except JAG.
\subjclass[2010]{05C25, 05C50, 20C30, 60K35}

%\date{}

%\dedicatory{}
%
%
%
\begin{abstract}
A conjecture by D. Aldous, which can be formulated as a statement
about the first nontrivial eigenvalue of the Laplacian of 
certain Cayley graphs on the symmetric group generated by transpositions,
has been recently proven by Caputo, Liggett and Richthammer.
Their proof is a subtle combination of two ingredients: 
a nonlinear mapping in the group algebra
of the symmetric group which permits a proof by induction, and
a quite hard estimate named the octopus inequality.
In this paper we present a simpler and more transparent proof
of the octopus inequality, which emerges naturally when looking
at the Aldous' conjecture from an algebraic perspective. 
We also show that the analogous of the Aldous' conjecture, 
where the spectral gap is replaced by the Kazhdan constant, 
does not hold in general.
\end{abstract}
%
% Keywords
% 05C25; Cayley graph; symmetric group; spectral gap; Aldous' conjecture; Kazhdan constant
%
\maketitle
\thispagestyle{empty}
%
%------------------------------------------------------------------
%
%
%
%\input sections.tex
\vspace{5mm}
%--------------------------------------------------------------------
\section{Introduction} 
%--------------------------------------------------------------------

\noindent
Let $G$ be a finite group with complex group algebra $\bC G$.
Given a representation $\dR$ of $G$ 
on the $d$-dimensional complex vector space $V$, and given an element 
$w = \sum_{g\in G} w_g \, g$ of the group algebra $\bC G$, 
we define the \textit{representation Laplacian} $\D_G(w,\dR)$
as the linear operator on $V$ given by
\begin{align}
  \label{eq:lapl}
  &\D_G(w,\dR) := \sum_{g\in G} w_g \, \[\id_V - \dR(g)\] = ( \dI_V - \dR)(w) 
  && w_g \in \bC
  ,
\end{align}
where $\id_V$ is the identity on $V$ and 
$\dI_V$ is the trivial representation given by $\dI_V(g) = \id_V$
for each $g\in G$. 
The \textit{support} of $w$ is defined as
\begin{equation*}
  \supp w := \{ g \in G : w_g \ne 0 \}
  \qquad
  w\in \bC G.
\end{equation*}
As a justification for using the term ``Laplacian'',
we observe that if $\dL$ is the left regular representation
of $G$ and $w = \sum_{q\in Q} q$, where
$Q\sset G$ is a symmetric generating set, then
$\D_G(w,\dL)$ is the standard (unnormalized) Laplacian of the Cayley graph
of $(G,Q)$. More generally, if $H$ is a subgroup of $G$, 
and $\dI \indR_H^G$ is the trivial
representation of $H$ induced to $G$, then
$\D_G( w, \dI \indR_H^G )$ is the Laplacian of the Schreier
graph of $(G,H,Q)$.

To the pair $(w,\dR)$ we can also associate a
\textit{spectral gap}, denoted by $\psi_G(w,\dR)$, which is essentially
the first
nontrivial eigenvalue of $\D_G(w,\dR)$ (see section 2 for a precise
definition). 
We then define the \textit{spectral gap of $w$},
by minimizing over representations
\begin{align}
  \label{eq:minrep}
  \psi_G(w) &:= \inf \bigl\{ \psi_G(w,\dR) : \dR \in \repf(G)
  \bigr\},
\end{align}
where $\repf(G)$ stands for the set of all the (equivalence classes
of) finite-dimensional representations of $G$ over the field of complex numbers.

Around 1992 a conjecture was formulated by David Aldous
about a problem which ``arose in conversation with Persi Diaconis''
\cite{Ald}. This conjecture, which became known as \textit{Aldous'
spectral gap conjecture}, originates in a probabilistic
framework and it asserts that two distinct Markov chains,
namely the random walk and the \textit{interchange process}
have generators with the same spectral gap. The reader is referred
to \cite{CaLiRi} for details on the probabilistic angle.
It is not too difficult to realize 
that this conjecture is actually a statement about the
representations of the symmetric group $S_n$
(see \cite{DiSh} or \cite{Ces1}). More precisely we have:

\begin{theorem}\label{thm:ald}
(Aldous' spectral gap conjecture, proven in \cite{CaLiRi}).
Let $S_n$ be the symmetric group, and let 
$w = \sum_{g\in S_n} w_g\, g \in \bC S_n$. 
Assume that 
\begin{enumerat}{(a)}
\item
$w_g\in \bR$ and $w_g\ge 0$ for each $g\in S_n$;
\item
$\supp w$ is a subset of
the set of all transpositions of $S_n$;
\item
$\supp w$ generates $S_n$.
\end{enumerat}
Then  we have
\begin{equation}
  \label{eq:D}
  \psi_G(w) = \psi_G(w, \dD_n) \,,
\end{equation}
where $\dD_n$ is the $n$-dimensional \textit{defining} representation
of $S_n$, associated with the natural action of
$S_n$ on the set $\{1,2,\ldots,n\}$.
\end{theorem}
Another way of stating this result (the spectral graph theory angle) 
is as follows:
let $\cG$ be a finite weighted graph
with vertex set $\{1,2,\ldots,n\}$ and symmetric weights $w_{ij} \ge 0$.
There is a natural way to associate a weighted Cayley graph 
$\Cay(\cG)$ to our original graph $\cG$:
any edge $e=\{i,j\}$ of $\cG$
can be identified
with a transposition $(ij)$ of the symmetric group $S_n$.
Consider then the Cayley graph with vertex set equal to $S_n$
and such that each edge $\{\pi, (ij) \pi \}$, where $\pi$ is a permutation
of $S_n$, carries a weight $w_{ij}$.
It is easy to show that the spectrum of the Laplacian of 
$\cG$ is a subset of the spectrum of the Laplacian of $\Cay(\cG)$.
Theorem \ref{thm:ald} is equivalent to the striking 
assertion that if $\cG$ is connected
(by edges with positive weight), then
the two Laplacian operators have
the same lowest nontrivial eigenvalue (see \cite{Ces1} for details).

Aldous' conjecture has been proven in several special cases by different authors 
in a series of papers \cite{DiSh}, \cite{FOW},  \cite{Bac}, \cite{HaJu},
\cite{KoNa}, \cite{Mor}, \cite{StCo}, \cite{Ces1} spanning about 25 years.%
\footnote{some of these papers even predate the ``official''
formulation of the conjecture, as it is often the case}
Unfortunately none of these papers contains tools which, 
with maybe some hammering and bending, can be 
forged to tackle the general case.
A proof which works in the general case has
required, in fact, a subtle combinations of two new ingredients: 
(a) a (nonlinear) ``reduction mapping'' 
$\th : \bC S_n \to \bC S_{n-1}$ which enables a proof
by induction on $n$; (b) an inequality, named the \textit{octopus inequality}
by the authors of \cite{CaLiRi}, 
who were ``inspired by its tentacular nature''.
Quite amazingly
both elements (a) and (b) appeared almost simultaneously
in the preprint versions of \cite{Diek} and \cite{CaLiRi}, but
the octopus inequality was left as a conjecture in \cite{Diek}
(it was actually proved in some particular cases), while it was proved
in \cite{CaLiRi}.

The main objective of this paper is to present (section \ref{sec:oct})
a simpler and more transparent
proof of the octopus inequality which hopefully can unravel
its tentacular nature. 
Our approach emerges naturally if one looks at paper \cite{CaLiRi}
with algebraically tinted glasses. For this reason we also 
include in section \ref{sec:ald} a more or less self-contained ``algebraic proof''
of Aldous' conjecture, showing how this conjecture follows
from the octopus inequality. The approach we follow in 
section \ref{sec:ald}, is not \textit{really} different from
what is done in \cite{CaLiRi}, since it uses the same reduction map $\th$. 
It is more like 
the same proof, as it appears from a different perspective.
We believe, however, that
this perspective may be useful for a better understanding
of what is going on and for investigating possible generalizations.

In the last section we discuss the connection 
with a result by M.~Kassabov \cite{Kas} where Aldous' conjecture
is proven%
\footnote{even if it is not explicitly mentioned.},
using a completely different approach,
for every finite Coxeter group $G$ in the special case of
$w = \sum_{q\in Q} q$, where $Q$ is a Coxeter generating set for $G$.
Kassabov also shows that, in such a case, the analogous 
of identity \eqref{eq:D} holds also for the Kazhdan constant. 
At the end of this section we also publicize an interesting
conjecture by P.~Caputo, which can be seen as a generalization
of Theorem \ref{thm:ald}.

\medno

%--------------------------------------------------------------------
\section{The representation Laplacian and its spectral gap} 
\label{sec:notation}
%--------------------------------------------------------------------

\noindent
Let $G$ be a finite group.
If $\dR$ is a finite-dimensional representation%
\footnote{representations are finite-dimensional throughout this
paper, with the exception of section \ref{sec:kazh}.}
of $G$ on the complex vector space $V$, we will sometimes
write that $(\dR, V)$ is a representation of $G$.
$\dR$  extends
to a representation of the group algebra $\bC G$ by letting,
for $w = \sum_{g\in G} w_g \, g$,
\begin{equation*}
  \dR(w) = \dR\Bigl( \sum_{g\in G} w_g \, g \Bigr) :=
  \sum_{g\in G} w_g \, \dR(g) 
  \qquad w_g\in\bC \,.
\end{equation*}
The degree of $\dR$ (that is the dimension of $V$) is denoted with $f_\dR$.
We denote with $\dI$ the
one-dimensional trivial representation, so that 
\begin{equation}
  \label{eq:triv}
  \dI(w) = \sum_{g\in G} w_g \in \bC \,.
\end{equation}
Since $G$ is finite, any representation $(\dR, V)$
can be realized as a \textit{unitary} representation.
%% This means that one can find an inner product on $V$
%% such that $\dR(g)$ is unitary w.r.t. such an inner product for each $g\in G$.
%% Viceversa if $\<\cdot,\cdot\>$
%% is a given inner product in $V$, then there is an equivalent
%% representation $\widetilde \dR$ which is unitary w.r.t. $\<\cdot, \cdot\>$.
$V^G$ stands for
the subspace of all invariant vectors
\begin{equation*}
  V^G := \{ v \in V : \dR(g) v = v,\ \forall g\in G \}.
\end{equation*}
We have that $V^G \ne \{0\}$ if and only if 
$\dI \sset \dR$, where, in general, $\dS \sset \dT$ means
that $\dS$ is a subrepresentation
of $\dR$. 
An eigenvalue $\l$ of the representation Laplacian 
$\D_G(w,\dR)$, defined in \eqref{eq:lapl}, will be called \textit{trivial}
if its corresponding eigenspace consists entirely
of invariant vectors $v \in V^G$. Clearly if $\l$ is trivial
then $\l=0$. The opposite implication will be discussed
in Proposition \ref{thm:lapl}.
We introduce a canonical involution in the group algebra $\bC G$ as
\begin{equation*}
  w = \sum_{g\in G} w_g \, g \too w^* := 
  \sum_{g\in G} \ol w_g \, g^{-1} 
\end{equation*}
and we denote the set of all symmetric elements as
\begin{equation*}
  \bC G^{(s)} := \{ w \in \bC G : w = w^* \} \,.
\end{equation*}
We will be also interested in a subset of $\bC^{(s)}$, \ie the
subset of all \textit{positive symmetric} elements, which
we denote with
\begin{equation*}
  \bC G^{(+)} := \{ w\in \bC G : w_g \in \bR,\ w_g \ge 0,\ 
  w_{g^{-1}} = w_g\ \forall g\in G \}.
\end{equation*}
In the following Proposition we summarize a few elementary properties
of the spectrum of the representation Laplacian.

\begin{proposition}\label{thm:lapl}
Let $\dR$ be a representation of $G$ on $V$, and
let $w\in \bC G$. We have
\begin{enumerat}{$(1)$}
\item
if $w \in \bC G^{(s)}$, then $\D_G(w,\dR)$ has real eigenvalues;
\item
if $w \in \bC G^{(+)}$, then $\D_G(w,\dR)$ has (real) non-negative 
eigenvalues;
\item
\label{it:3}
let $w \in \bC G^{(+)}$ and assume that $\supp w$ generates $G$. 
If $0$ is an eigenvalue of $\D_G(w,\dR)$, then its
eigenspace  coincides with $V^G$,
hence it is necessarily trivial.
\end{enumerat}
\end{proposition}

\Pro\ 
Let $\<\cdot,\cdot\>$ be an inner product on $V$ which makes
$\dR$ a unitary representation.
If $w\in \bC G^{(s)}$, 
then we have
\begin{align*}
  \dR(w) = 
  \dR(w^*) = \sum_{g\in G} \ol w_g \, \dR( g^{-1} ) 
  = \sum_{g\in G} \ol w_g \, \dR( g )^* = \dR(w)^* \,, 
\end{align*}
\ie $\dR(w)$ is a self-adjoint linear operator, which proves (1).

\smallno
Since $\dR$ is unitary and the coefficients $w_g$
are real and non-negative, from \eqref{eq:lapl} we obtain for $v\in V$
\begin{equation}
\begin{split}
  \label{eq:pos}
  \< \D_G(w,\dR)\,  v, v \> &=
  \sum_{g\in G} w_g \[ \|v\|^2 - \< \dR(g) v,v \> \]
  \\ &
  =
  \ov2 \sum_{g\in G} w_g \[ 2 \|v\|^2 - \< (\dR(g) + \dR(g^{-1}) \, v,v \> \]
  \\ &
  =
  \ov2 \sum_{g\in G} w_g \[ 2 \|v\|^2 - \< \dR(g) v,v \>  - \< v, \dR(g) v \> \]
  \\ &
  =
  \ov2 \sum_{g\in G} w_g \, \| \dR(g) v - v \|^2
  \ge 0 .
\end{split}
\end{equation}
Hence $\D_G(w,\dR)$ is positive semidefinite as quadratic form, which
implies (2).

\smallno
For proving (3) let $v \ne 0$ with
$\D_G(w,\dR) v = 0$. Then the last inequality in
\eqref{eq:pos} is in fact an equality. By consequence
$\dR(g) v = v$ for each
$g\in \supp w$. Since $\supp w$ generates $G$, we 
conclude that $v\in V^G$. 

The 
inverse implication, that is if $v$ is an invariant
vector then $v$ belongs to the $0$-eigenspace of $\D_G(w,\dR)$
is trivial.
\qed

\medno
%% In general $\D(w,\dR)$ may have a 
%% nontrivial zero eigenvalue, but if $\supp w$ generates $G$
%% this does not happen.
If $w \in \bC G^{(+)}$,
we define the \textit{spectral gap}
of the pair $(w, \dR)$ as
\begin{align}
  \label{eq:sgdef}
  \psi_G(w, \dR ) &:= \min \{ \l \in \spec \D_G(w,\dR) :
  \text{$\l$ is nontrivial} \}
  && w \in \bC G^{(+)} 
  ,
\end{align}
with the convention that $\min \emp = +\oo$.%
\footnote{This happens when $\dR$ is a multiple of $\dI$.
In this case $0$ is the only eigenvalue and it is trivial.}
The \textit{spectral gap of $w$} is defined by minimizing
over representations
\begin{align}
  \label{eq:minrep1}
  \psi_G(w) &:= \inf \bigl\{ \psi_G(w,\dR) : 
  \dR \in \repf(G)
  \bigr\} \,.
\end{align}
By Maschke's complete reducibility theorem, 
we have for each $\dR \in \repf(G)$,
\begin{equation}\label{eq:Y}
  \dR \cong \bigoplus_{\dT \in \Irr(G)} r(\dT) \, \dT \,,
\end{equation}
where $\Irr(G)$ stands for the set of all the (equivalence classes
of) irreducible complex representations of $G$,
and $r(\dT)$ are suitable nonnegative integers.
By consequence
\begin{align}\label{eq:Y1}
  \spec \D_G( w,\dR) = \bigcup_{\dT\in \Irr(G) :\, r(\dT)>0} 
  \spec \D_G( w,\dT ) \,,
\end{align}
which implies
\begin{equation*}
  \psi_G(w,\dR) = \min \{ 
  \psi_G(w, \dT ) 
  : \dT\in \Irr(G) :\, r(\dT)>0 \}.
\end{equation*}
For this reason in \eqref{eq:minrep1} we can limit ourselves
to consider irreducible representations, that is
\begin{align}
  \label{eq:minrep2}
  \psi_G(w) &= \min \bigl\{ \psi_G(w,\dR) : 
  \dR \in \Irr(G)
  \bigr\} \,,
\end{align}
where we have also replaced $\inf$ with $\min$ since
$\Irr(G)$ has finite cardinality.
Furthermore, because of our definition \eqref{eq:sgdef}, we have
$\psi_G(w,\dI) = +\oo$, hence the conditions $\dR\ne \dI$ and $\dI \not\sset \dR$
be freely added in \eqref{eq:minrep2} and 
in \eqref{eq:minrep1} respectively.
Let $\dL$ be the left regular representation of $G$.
Since $\dL$ contains all irreducible representations, 
\begin{equation}
  \label{eq:L}
  \dL = \bigoplus_{\dT \in \Irr(G)} f_\dT \, \dT
\end{equation}
$\dL$
determines the spectral gap of $w$, that is
$\psi_G(w) = \psi_G(w,\dL)$.
Statement \ref{it:3} of Proposition \ref{thm:lapl} implies that
if $w\in \bC G^{(+)}$ and $\supp w$ generates $G$, then
$\psi_G(w) > 0$.
If $A$ is a self-adjoint linear operator on a finite dimensional
vector space we also define
\begin{equation*}
  \l^*(A) := \max \spec A \,.
\end{equation*}
If $\dI \not\subset \dR$, then
$\dR$ has no invariant nonzero vectors, thus $\D_G(w,\dR)$ has no trivial
eigenvalue, which, in view of \eqref{eq:lapl}, implies
\begin{equation}
  \label{eq:IT}
  \psi_G(w,\dR)  = \min \spec \D_G(w,\dR) 
  = \dI(w) - \l^*( \dR(w) ) \,.
\end{equation}

\medno

%--------------------------------------------------------------------
\section{Aldous' conjecture: an algebraic perspective} 
\label{sec:ald}
%--------------------------------------------------------------------

\noindent
In this section we prove Theorem \ref{thm:ald}, given the
octopus inequality whose proof is postponed to section \ref{sec:oct}.
The fundamental tool is a map $\th : \bC S_n \to \bC S_{n-1}$ 
introduced in \cite{CaLiRi} and \cite{Diek}.

For a self-adjoint linear operator $A$ we write $A\ge 0$ if 
$\<A \cdot, \cdot\>$ is
a positive semidefinite bilinear form.
For a finite group $G$,
we let
\begin{equation}
\begin{split}
  \label{eq:G}
  \G(G) := 
  &\{ w \in \bC G^{(s)} : 
  \D_{G}(w, \dT) \ge 0, 
  \ \forall \dT\in \Irr(G)  
  \} \\
  =
  &\{ w \in \bC G^{(s)} : 
  \D_{G}(w, \dL) \ge 0
  \} 
  \,,
\end{split}
\end{equation}
where, as usual, $\dL$ denotes the left regular representation.
The two definitions are equivalent thanks to \eqref{eq:L}.
We know the $\bC G^{(+)}$ is a subset of $\G(G)$, but we will have
to deal with elements $w$ with negative coefficients.
A necessary and sufficient condition for $w \in \G(G)$ is, clearly
\begin{equation}
  \label{eq:l*}
  \l^*( \dT(w) ) \le \dI(w)
  \qquad
  \forall \dT\in \Irr(G) \,.
\end{equation}
If $H$ is a subgroup of $G$ we denote with
$\dT \resR^G_H$ the restriction of the representation $\dT$
to $H$.

$\G(G)$ plays a central role in the strategy for estimating
the spectral gap. The idea can be roughly
described as follows: supposed we are interested in a lower bound
on $\psi_G(w)$. Let then $H$ be a subgroup of $G$ and let 
$z \in \bC H^{(+)}$. If $w-z \in \G(G)$, then it turns out
that $\psi_G(w)$ can ``almost'' be estimated in terms of
$\psi_H(z)$, where ``almost'' means that one has still to
worry about irreducible representations of $G$ which, 
when restricted to $H$, contain the trivial representation.
This idea is implemented in the next two results.

\begin{lemma}\label{thm:gg1}
Let $G$ be a finite group and let $H$ be a subgroup of $G$.
If $w\in \bC G^{(+)}$ and $z\in \bC H^{(+)}$ are such that
$w-z \in \G( G )$, then, for any irreducible representation
$\dT$ of $G$ such that $\dI \not\subset \dT\resR^G_H$, we have
\begin{align}
  \label{eq:psi}
  \psi_G(w, \dT) \ge 
  \min \bigl\{ 
    \psi_H( z, \dS) : \dS\in \Irr(H),\ \dS \sset \dT \resR^G_H 
  \bigr\}
  \,.
\end{align}
\end{lemma}

\Pro\ 
Let $y := w-z \in \G(G)$, and let $\dT$ be an irreducible
representation of $G$.
Since $z \in \bC H^{(s)}$ we have
\begin{equation*}
  \l^*( \dT( z ) ) = \l^*\Bigl( \dT \resR^G_H(z) \Bigr)
  =  \max\bigl\{ \l^*( \dS(z) ) :  \dS\in \Irr(H),\ \dS \sset \dT \resR^G_H  \bigr\}
  \,.
\end{equation*}
Since $\l^*(A+B) \le \l^*(A) + \l^*(B)$
for each pair of self-adjoint linear operators $A,B$,
and using the fact that
$y\in \G(G)$, so that \eqref{eq:l*} applies to $y$, we get
\begin{equation}
\begin{split}
  \label{eq:gg1}
  \l^*( \dT(w) ) &= 
  \l^*( \dT(z) + \dT(y) ) 
  \le
  \l^*( \dT(z) ) + \l^*( \dT(y) ) 
  \\ &
  \le 
  \l^*( \dT(z) ) + \dI(y)
  \\ &
  = \max \bigl\{ \l^*( \dS(z)) : 
  \dS\in \Irr(H),\ \dS \subset \dT \resR^G_H \bigr\} 
  + \dI(w) - \dI(z)
  \,.
\end{split}\end{equation}
Since, by hypothesis, no subrepresentation $\dS$ 
of $\dT \resR^G_H$ is equal to the trivial one,
we have that \eqref{eq:psi} follows from \eqref{eq:gg1} and \eqref{eq:IT}. \qed

\medno
As a consequence we get the following ``semi-recursive'' result:

\begin{proposition}\label{thm:sr}
Let $w$ and $z$ be as in Lemma \ref{thm:gg1}. Then
\begin{align*}
  \psi_G(w) \ge \min\Bigl\{ \psi_H(z),\, 
  \min \bigl\{ \psi_G( w,\dT) : \dT \in \Irr(G),\ 
  \dT\resR^G_H \supset \dI \bigr\}\; \Bigr\} 
  \,.
\end{align*}
\end{proposition}

\Pro\ 
Let us write the set $\Irr(G)$ as a disjoint union
$\Irr(G) = \cA \cup \cB$,
where
\begin{align*}
  \cA &:= \{ \dT \in \Irr(G) : \dT\resR^G_H \supset \dI \}
  \\
  \cB &:= \{ \dT \in \Irr(G) : \dT\resR^G_H \not\supset \dI \}
  \,.
\end{align*}
We can now apply Lemma \ref{thm:gg1} to each $\psi_G(w,\dT)$ 
when $\dT\in \cB$. In this way we obtain
\begin{align*}
  \min_{\dT\in \cB} \psi_G(w,\dT) \ge
  \min_{ \dS \in \Irr(H),\ \dS\ne \dI } \psi_H(z,\dS) =: \psi_H(z)
  \,,
\end{align*}
and the Proposition follows.
\qed

\medno
\textit{The case of the symmetric group}.
Let now consider the special case $G=S_n$.
We write $\a\partit n$ if $\a = (\a_1, \a_2 \, \ldots)$ 
is a (proper) partition of $n$. The equivalence classes of 
irreducible representations
of $S_n$ are in one-to-one correspondence with the set of all
partitions of $n$. We denote with $[\a]$ the equivalence class
corresponding to partition $\a$. 
The defining representation $\dD_n$ of the symmetric group $S_n$
decomposes as 
\begin{equation}
    \label{eq:DI}
    \dD_n = \dI \oplus [n-1,1]
\end{equation}
as a direct sum of irreducible
factors, hence 
\begin{equation*}
  \psi_{S_n}(w, \dD_n) = \psi_{S_n}(w, [n-1,1]).
\end{equation*}
If $G=S_n$ and $H=S_{n-1}$, the only irreducible 
representations of $S_n$ which contain the trivial representation when
restricted to $S_{n-1}$ are $\dI$ 
and $[n-1,1]$ (see \eqref{eq:br}). 
Hence, in this particular case, Proposition \ref{thm:sr} becomes

\begin{proposition}\label{thm:sr1}
Let $w\in \bC S_n^{(+)}$, $z\in \bC S_{n-1}^{(+)}$ be such that $w-z \in \G(S_n)$.
Then
\begin{align*}
  \psi_{S_n}(w) \ge \min\bigl\{ \psi_{S_{n-1}}(z),\, 
  \psi_{S_n}( w, \dD_n ) 
  \bigr\} 
  \,.
\end{align*}
\end{proposition}

\medno

\subsection{Proof of Aldous' conjecture}

Imagine that we want to prove \eqref{eq:D} for all $w$
in a certain subset $\cA_n \sset \bC S_n^{(+)}$, \textit{for all $n$}.
The strategy is as follows: 
suppose
that, for all $n\ge 3$,
we find a map $\th : \cA_n \to \cA_{n-1}$ such that,
for all $w\in \cA_n$, we have
\begin{enumerate}[(a)]
\item $w - \th(w) \in \G(S_n)$ 
\item 
%% REFERENCE
$\psi_{S_n}( w, \dD_n ) \le \psi_{S_{n-1}}(  \th(w), \dD_{n-1})$.
\end{enumerate}
Then the proof can be done by induction
(assuming that we have a good starting point at $n=2$): in fact, if
\eqref{eq:D} holds for $n=k-1$, then
\begin{align}
    \label{eq:ind}
  &\psi_{S_{k-1}}(z) = \psi_{S_{k-1}}( z, \dD_{k-1} )
  &&\forall z\in \cA_{k-1} \,.
\end{align}
From Proposition \ref{thm:sr1} and \eqref{eq:ind} with $z=\th(w)$, 
and from property (b) of the map $\th$
it follows that, for all $w\in \cA_k$, we have
\begin{align*}
  \psi_{S_k}(w) &\ge \min\bigl\{ 
  \psi_{S_{k-1}}( \th(w), \dD_{k-1}  ) \,, \; 
  \psi_{S_k}( w, \dD_k ) 
  \bigr\}
  \\
  &=
  \psi_{S_k}( w, \dD_k ) \,,
\end{align*}
which implies $\psi_{S_k}(w) = \psi_{S_k}( w, \dD_k )$. 
Hence the induction is completed.

\smallno
Let now $\cA_n$ be the set considered in Theorem \ref{thm:ald},
\ie the set of all $w \in \bC T_n^{(+)}$ such that
$\supp w$ generates $S_n$. 
The starting point of the induction is trivial because
the only irreducible representations of $S_2$ are $[2]$ and $[1,1]$.
The map $\th$ which does the trick
was found in \cite{CaLiRi} and \cite{Diek} and is given by
\begin{equation*}
  \th : w = \sum_{(ik)\in T_n} w_{ik}(ik) \longmapsto  
  \sum_{ (ik) \in T_{n-1}} 
  \[ w_{ik} + \frac{w_{in} \, w_{kn}}{w_{1n} + \cdots + w_{n-1,n}} \] (ik)
  \,.
\end{equation*}
All is left is proving properties (a) and (b) for this map.
Let, for brevity,  
\begin{equation*}
  \D_n := \D_{S_n}(w, \dD_n) \qquad  \D_n^\th := \D_{S_n}( \th(w), \dD_n) \,.
\end{equation*}
It is straightforward to check that the matrix elements of $\D_n - \D_n^\th$ are given by
\begin{equation*}
  \bigl[ \D_n - \D_n^\th \bigr]_{ij} = \frac{\d_i \,\d_j}{\d_n} \,,
\end{equation*}
where $\d_i = -w_{in}$ for $i=1,\ldots,n-1$ and $\d_n = \sum_{i=1}^{n-1} w_{in}$.
Following \cite{Diek} we realize that $\D_n-\D_n^\th$
is a positive semidefinite rank-1 matrix, so by standard linear algebra results
as \cite[Cor 4.3.3 and Thm. 4.3.4]{HorJohn}, one obtains the following
interlacing property of the eigenvalues: 
let $\l_{n,k}$ and $\l_{n,k}^\th$ denote the $k$-th lowest eigenvalue of
respectively $\D_n$ and $\D_n^\th$. Then
\begin{equation}
  \label{eq:inter}
  \l^\th_{n,1} \le \l_{n,1} \le
  \l^\th_{n,2} \le \l_{n,2} \le
  \cdots \le
  \l^\th_{n,n} \le \l_{n,n} \,.
\end{equation}
Since $\th(w) \in \bC S_{n-1}$, the last row and 
the last column of $\D_n^\th$ are zero, more precisely
we have 
\begin{equation*}
  \D_n^\th =   \D_{n-1}^\th \oplus [0]_{1\times 1} 
\end{equation*}
where $[0]_{1\times 1}$ is the $1\times 1$ zero matrix. 
Thus
\begin{equation*}
  \l_{n,k}^\th = \l_{n-1, k-1}^\th
  \qquad
  k=1, \ldots, n\,.
\end{equation*}
This implies that the spectral gap of $\D_{n-1}^\th$
corresponds to the third lowest eigenvalue of $\D_n^\th$.
Combining with \eqref{eq:inter} we get 
\begin{equation*}
  \psi_{S_n}(w, \dD_n) 
  = \l_{n,2} \le \l_{n,3}^\th = \l_{n-1,2}^\th = 
  \psi_{S_{n-1}}(\th(w), \dD_{n-1}) \,,
\end{equation*}
thus property (b) holds. 

The proof of property (a) is definitely harder and it corresponds
to proving the octopus inequality which we do in the next section,
Theorem \ref{thm:wg}.
\qed

\medno

%--------------------------------------------------------------------
\section{The Octopus} 
\label{sec:oct}
%--------------------------------------------------------------------

\noindent
In this section we prove property (a) of the map $\th$ defined in
the previous section, \ie we show that for each $n\ge 3$ we have
\begin{equation*}
  \D_{S_n}( w - \th(w), \dL) \ge 0
  \qquad
  \forall w \in \bC T_n^{(+)} \,.
\end{equation*}
It follows from the definition of the left regular representation $\dL$
that this is equivalent to the octopus inequality as
stated in \cite[Thm. 2.3]{CaLiRi} or \cite[Conjecture 1]{Diek}.

We start by recalling a few well known facts about $S_n$ that we need later.
The symmetric group $S_n$ is the (disjoint) union of its conjugacy classes
\begin{equation*}
  S_n = \bigcup_{\a\partit n} C^\a \,,
\end{equation*}
where, for $\a = (\a_1, \a_2, \ldots) \partit n$,  $C^\a$ is the set of all
permutations $\pi$ which are product of disjoint cycles of lengths 
$\a_1$, $\a_2$, \ldots.
The partition $\a$ is called the \textit{cycle partition} of $\pi$.
The cardinality of these conjugacy classes is given \cite[Ch.1]{Sag} by
\begin{align}
  \label{eq:conj}
  |C^\a| = n! \; 
  \Bigl[ \prod_{k=1,2,\ldots} k^{\a_k^\#} \, (\a_k^\#) ! \; \Bigr]^{-1}
\end{align}
where we denote with $\a^\#$ the \textit{content} of $\a$, given by
\begin{equation*}
   \a^\#_k := | \{ i : \a_i = k \} | 
  \qquad k=1,2,\ldots \,.
\end{equation*}
For instance, if $\a = (3,3,2,1,1,1)$, we have $\a^\# = (3,1,2)$ and
\begin{equation*}
  |C^\a| = \frac{ 11! }{ 1^3\, 3!\, 2^1\, 1!\, 3^2\, 2! } \,.
\end{equation*}

\medno
An irreducible representation $[\a]$ of $S_n$ is no longer
(in general) irreducible when restricted to a subgroup $H$
of $S_n$. If $H = S_{n-1}$ there is a simple
\textit{branching rule} \cite[section 2.8]{Sag}
for computing the coefficient of the decomposition:
\begin{align}
  \label{eq:br}
  [\a] \resR^{S_n}_{S_{n-1}} = \bigoplus_{\b \in \a^-} \, [\b]
  \qquad \a\partit n
\end{align}
where, if $\a=(\a_1, \ldots,\a_r)$, 
 $\a^-$ is defined as the collection of all sequences
of the form
\begin{equation*}
  (\a_1, \ldots, \a_{i-1}, \a_i -1, \a_{i+1}, \ldots, \a_r)
\end{equation*}
\textit{which are partitions of} $n-1$. 
For example
\begin{align*}
  [6,5,5,3,1]\resR^{S_{20}}_{S_{19}} = 
  [5,5,5,3,1]
  \oplus [6,5,4,3,1]
  \oplus [6,5,5,2,1]
  \oplus [6,5,5,3] \,.
\end{align*}

\medno
For $\a\partit n$, we define the following element in the group algebra of $S_n$
\begin{align}
  \label{eq:Ja}
  \bC S_n \ni J^\a := \sum_{\pi \in C^\a}  \pi \,.
\end{align}
If $A$ is a subset of $\{1,\ldots,n\}$ with $|A|=m$, and $\a\partit m$, we also let
\begin{equation}
  \label{eq:Jaa}
  J^\a_A := \sum_{\pi \in C^\a_A} \pi \,,
\end{equation}
where $C^\a_A$ is the set of all permutations $\pi \in S_A$ whose cycle partition
is equal to $\a$, and
$S_A$ stands for the group of permutations
of $A$. Since $\bC S_A$ is naturally embedded in $\bC S_n$, we can regard
$J^\a_A$ as an element of $\bC S_n$ if useful.

For $\a,\b \partit n$, we denote with $\c^\b(\a)$ the value taken by the
irreducible character $\c^\b$ on the conjugacy class $C^\a$, while $f_\b$
stands for the degree of the irreducible representation $\b$.

A more or less standard application of Schur's Lemma implies that
the image of $J^\a$ under any irreducible representation of $S_n$
is a multiple of the identity. More precisely we have:

\begin{lemma}\label{thm:schur}
If $\a, \b$ are two partitions of $n$, then
\begin{equation*}
  \dT^\b( J^\a ) = \frac{ |C^\a| \, \c^\b(\a) }{ f_\b } I_{ f_\b } \,.
\end{equation*}
\end{lemma}

\Pro\ 
Since $\pi J^\a = J^\a \pi$ for each $\pi\in S_n$, by Schur's Lemma
$\dT^\b( J^\a ) = c I_{f_\b}$ for some $c\in \bC$. Taking the trace
of both sides produces the result. \qed

\begin{theorem}\label{thm:wg}
For each $w \in \bC T_n^{(+)}$ we have $w - \th (w) \in \G(S_n)$.
\end{theorem}

\Pro\ 
We start with the elementary observation that
\begin{equation}
\begin{split}
  \label{eq:clo}
  &\text{$\G(S_n)$ is closed under linear combinations} \\[-1mm]
  &\text{with nonnegative coefficients.}
\end{split}
\end{equation}
Let $w\in \bC T_n^{(+)}$, thus
\begin{equation*}
  w = \sum_{ (ij) \in T_n } w_{ij} \, (ij)
  \qquad w_{ij} \ge 0 \,.
\end{equation*}
Let $\hw := \dI(w) \, (w - \th w)$. Letting, for simplicity, $x_i := w_{in}$
for $i=1, \ldots, n-1$, we have
\begin{equation}
\begin{split}
  \label{eq:hw}
  \hw &= \( \sum_{i=1}^{n-1} x_i \) \sum_{j=1}^{n-1} x_j\,(jn)
  - \ov2 \sum_{ \atopp{i,j=1}{i\ne j} }^{n-1} x_i x_j \, (ij)
  \\
  &=
  \sum_{i=1}^{n-1} x_i^2 \, (in) + \sum_{ \atopp{i,j=1}{i < j} }^{n-1} 
  x_i x_j \,[ (in) +(jn) -(ij) ] \,.
\end{split}
\end{equation}
The idea of the proof is as follows:
$\hw^2$ is a homogeneous quartic polynomial in the variables $(x_i)_{i=1}^{n-1}$
with coefficients in the group algebra $\bC S_n$
\begin{equation*}
  \hw^2 = 
  \sum_{ \mu \in \bN^{n-1} :\, |\mu| = 4  } B_\mu^{(n)}  x^\mu 
  \qquad B^{(n)}_\mu \in \bC S_n \,,
\end{equation*}
where $\mu$ is a multi-index, $|\mu| = \sum_i \mu_i$ and $x^\mu = \prod_{i} x_i^{\mu_i}$.
If we can show that
\begin{equation}
  \label{eq:Bn}
  \text{ for each $\mu\in \bN^{n-1}$ with $|\mu|=4$ we have $B^{(n)}_\mu \in \G(S_n)$ }
\end{equation}
then the theorem is proved. 
In fact, thanks to \eqref{eq:clo}, 
since $x^\mu \ge 0$, 
statement \eqref{eq:Bn} implies
that $\hw^2\in \G(S_n)$. Thus, by \eqref{eq:l*}, we obtain
\begin{equation*}
  \dT^\a(\hw)^2 = \dT^\a(\hw^2) \le \dI( \hw^2 ) \, I_{f_\a} = \dI( \hw )^2 \, I_{f_\a}  
  \qquad \forall \a \partit n \,.
\end{equation*}
Since $\dI(w) > 0$, it follows that each eigenvalue $\l$ of $\dT^\a(\hw)$ 
satisfies $|\l| \le \dI(\hw)$, hence $\hw \in \G(S_n)$. 
Since $w-\th(w) = c\, \hw$ with $c>0$, we conclude that $w-\th(w) \in \G(S_n)$.

\begin{remark}\label{thm:compare}
A somehow similar approach is pursued in
\cite[Lemma 3.1, 3.2]{CaLiRi}, where it is proved that $\D_{S_n}(\hw^2, \dL) = -
  \sum_{ \mu \in \bN^{n} :\, |\mu| = 4  } A_\mu^{(n)}  x^\mu$,
where $A_\mu^{(n)}$ are positive semidefinite matrices.
The difference is that this expansion is with respect to
the $n$ \textit{dependent} variables  $x_0, x_1, \ldots, x_{n-1}$, where
$x_0 := - \sum_{i=1}^{n-1} x_i$. By consequence the monomials
$x^\mu$ have no definite sign and thus one cannot simply get rid of them.
The required extra work is done in Lemma 3.3 and Lemma 3.4 of \cite{CaLiRi}.
\end{remark}

\medno
We are thus left with the proof of \eqref{eq:Bn}. The first step
is to find manageable expressions for the coefficients $B^{(n)}_\mu$.
The next Lemma does not require the coefficients $x_i$ to be non-negative.

\begin{lemma}
\label{thm:w2}
If $w\in \bC T_n$ and $\hw = \dI(w) (w-\th w)$, then we have%
\footnote{if $n<5$ one or more of the following terms are missing, which is ok.}
\begin{equation}
\begin{split}
  \hw^2 &= \Biggl( \, \sum_{i=1}^{n-1} x_i^4 
  + 2 \sum_{ \atopp{i,j=1}{i \ne j} }^{n-1} x_i^3 x_j 
  + 3 \sum_{ \atopp{i,j=1}{i < j} }^{n-1} x_i^2 x_j^2 \Biggr) \cdot \uuno_{S_n} 
  \\ &
  + \!\!\! \sum_{ \atopp{i,j,k=1}{i\ne j,\, i \ne k,\, j<k} }^{n-1} \!\!\! x_i^2 x_j x_k
  \Bigl( X_{\{i,j,k,n\}} + 2 \cdot \uuno_{S_n} \Bigr) 
  + 2 \!\sum_{ \atopp{i,j,k,\ell=1}{i <j<k<\ell} }^{n-1} x_i x_j x_k x_\ell \, Y^{n}_{\{i,j,k,\ell\}}
\end{split}
\end{equation}
where $\uuno_{S_n}$ is the identity in $S_n$ and, remembering \eqref{eq:Jaa},
\begin{align}
  \label{eq:X}
  X_{\{i,j,k,n\}} &:= J^{(3,1)}_{\{i,j,k,n\}} - 2 J^{(2,2)}_{\{i,j,k,n\}} \\
  \label{eq:hY}
  Y^n_{\{i,j,k,\ell\}} &:= J^{(3,1,1)}_{\{i,j,k,\ell, n\}} - J^{(2,2,1)}_{\{i,j,k,\ell, n\}} 
  - X_{\{i,j,k,\ell\}} \,.
\end{align}
\end{lemma}

\medno
We postpone the proof of this lemma and we observe that, thanks to \eqref{eq:clo},
in order to show that $\hhW^2 \in \G(S_n)$, and thus conclude the proof
of Theorem \ref{thm:wg}, it is sufficient to show that both $X_{\{i,j,k,n\}}$ 
and $Y^n_{\{i,j,k,\ell\}}$
are in $\G(S_n)$. We also claim that it is sufficient to consider two special cases
\begin{equation}
  \label{eq:n4-5}
  X_{\{1,2,3,4\}} \in \G(S_4) \qquad \text{and}\qquad Y^5_{\{1,2,3,4\}} \in \G(S_5)\,.
\end{equation}
Let in fact $\s, \pi$ be 2 permutations of $S_n$. If $\s$ is written as a product
of disjoint cycles, then $\pi \s \pi^{-1}$ is the permutation which is obtained from $\s$
by replacing in each cycle $i$ with $\pi(i)$. Thus, if $\pi$ sends $(i,j,k,n)$ 
to $(1,2,3,4)$, then 
\begin{equation*}
  \pi \, X_{\{i,j,k,n\}} \, \pi^{-1} = X_{\{1,2,3,4\}}.
\end{equation*} 
By consequence, for any $\dR \in \repf(S_n)$, the two matrices
$\dR( X_{\{i,j,k,n\}} )$ and $\dR( X_{\{1,2,3,4\}} ) $ are equivalent.
We then observe that, in order to prove that $X_{\{1,2,3,4\}} \in \G(S_n)$,
it is sufficient to consider the irreducible representations of $S_4$. 
Let in fact $\a\partit n$ and let $\dT^\a$ be the corresponding irreducible representation
of $S_n$. 
Then we have
\begin{equation*}
  \dT^\a\resR^{S_n}_{S_4} = \bigoplus_{\b \partit 4} \nu_\b \, \dT^b
\end{equation*}
for some
nonnegative integral coefficients $\nu_\b$, 
Since $X_{\{1,2,3,4\}}$ is an element of $S_4$, we get
\begin{equation*}
  \dT^\a( X_{\{1,2,3,4\}} ) = \bigoplus_{\b \partit 4} \nu_\b \, \dT^\b(  X_{\{1,2,3,4\}} ) \,,
\end{equation*}
and so the spectrum of $ \dT^\a( X_{\{1,2,3,4\}} )$ is the union of the spectra
of those $ \dT^\b( X_{\{1,2,3,4\}} )$, where $\b$ is a partition of $4$ 
for which $\nu_\b > 0$. This shows that if $X_{\{1,2,3,4\}}$ belongs to $\G(S_4)$
it necessarily belongs to $\G(S_n)$ for all $n \ge 4$.

\smallno
A similar argument applies to $Y^{n}_{\{i,j,k,\ell\}}$. 
In the following we let for simplicity
\begin{align*}
  &\hX := X_{\{1,2,3,4\}} 
  &&\hY := Y^5_{\{1,2,3,4\}}.
\end{align*}
So the theorem is proved is we show that \eqref{eq:n4-5} hold, \ie 
if we prove that
\begin{align}
    \label{eq:n4}
  \l^*\bigl[ \dT^\a( \hX ) \bigr] &\le \dI( \hX ) &&\forall \a\partit 4 \\
    \label{eq:n5}
  \l^*\bigl[ \dT^\a( \hY ) \bigr] &\le \dI( \hY ) &&\forall \a\partit 5 \,.
\end{align}
It is well known 
\cite[2.1.8]{JaKe} that if $z\in \bC S_n$ is a linear combinations
of \textit{even} permutations and  $\a'$ if the conjugate  partition%
\footnote{which means that the Young diagram of $\a'$ is obtained
from the Young diagram of $\a$ by interchanging rows and columns.} 
of $\a$, 
then $\dT^{\a'}(z) \cong \dT^{\a}(z)$, so we must only check previous inequalities
for (roughly) half the irreducible representations of $S_4$ and $S_5$.
Thanks to Lemma \ref{thm:schur} and to \eqref{eq:conj} we have
\begin{align*}
  \dT^\a( \hX ) &= 
  f_\a^{-1} \[
  |C^{(3,1)}| \, \c^\a(3,1) 
  - 2 \, |C^{(2,2)}| \, \c^\a(2,2) \]  \cdot I_{f_\a}
  \\ &
  = f_\a^{-1} \[
  8 \, \c^\a(3,1) 
  - 6 \, \c^\a(2,2) \] \cdot I_{f_\a} \,.
\end{align*}
Hence $\dT^\a(\hX)$ is a multiple of the identity, and, in particular, 
\begin{align}
  \label{eq:lX}
  X^\a := \l^*( \dT^\a(\hX) ) = f_\a^{-1} \[
  8 \, \c^\a(3,1) 
  - 6 \, \c^\a(2,2) \] \,.
\end{align}
The relevant characters can be quickly computed 
with the Murnaghan-Nakayama rule.%
\footnote{or, even more quickly, one can look them up on the internet.}
Moreover
\begin{align*}
  \dI( \hX ) = |C^{(3,1)}| - 2\, |C^{(2,2)}| = 2 \,.
\end{align*}
Table \ref{tab:n4} shows that 
\eqref{eq:n4} holds.%
\footnote{Partitions $\a=(2,1,1)$ and $\a=(1^4)$, absent in this table, are
conjugate to $(3,1)$ and $(4)$ respectively.}

\smallno
For $\hY$ things are slightly more complicated.
If $\a\partit 5$, since $\hX\in S_4$, we have
\begin{equation*}
  \dT^\a(\hX) = \dT^\a \resR^{S_5}_{S_4}(\hX) = \bigoplus_{\b\in \a^-} \dT^\b( \hX ) \,.
\end{equation*}
Since $|C^{(3,1,1)}| = 20$ and $|C^{(2,2,1)}| = 15$,
from \eqref{eq:hY} and \eqref{eq:br} we obtain
\begin{align*}
  \dT^\a( \hY ) &= 
  f_\a^{-1} \[ 20 \, \c^\a(3,1,1) 
  - 15 \, \c^\a(2,2,1) \]  \cdot I_{f_\a} - \oplus_{\b\in \a^-} \dT^\b( \hX ) \,,
\end{align*}
which, thanks to \eqref{eq:lX}, implies
\begin{align*}
  Y^\a := \l^* \bigl[\dT^\a( \hY ) \bigr] &= 
  \max_{\b \in \a^-} F_{\a \b} \,,
\end{align*}
where
\begin{align*}
    F_{\a \b} &:=
  f_\a^{-1} \[ 20 \, \c^\a(3,1,1) 
  - 15 \, \c^\a(2,2,1) \] - X^\b \,.
\end{align*}
On the other hand, we have
\begin{align*}
  \dI( \hY ) = |C^{(3,1,1)} | - |C^{(2,2,1)}| - \dI( \hX ) = 3 \,.
\end{align*}
We collect all relevant quantities in table \ref{tab:n5} which shows that
\eqref{eq:n5} holds.

\begin{table}
\begin{center}
\caption{Proof of \eqref{eq:n4}}
\label{tab:n4}%
\renewcommand\arraystretch{1.2}
\begin{tabular}{|c|c|c|c|c|} \hline 
$\a$    & $f_\a$  & $\c^\a(3,1)$  & $\c^\a(2,2)$  & $ X^\a $ 
\\ \hline \hline
(4)     &  1      &  1            &  1            &  2      \\ \hline
(3,1)   &  3      &  0            & -1            &  2      \\ \hline
(2,2)   &  2      & -1            &  2            &-10      \\ \hline
\end{tabular}
\end{center}
\end{table}

\begin{table}
\begin{center}
\caption{Proof of \eqref{eq:n5}}
\label{tab:n5}%
\renewcommand\arraystretch{1.2}
\begin{tabular}{|c|c|c|c|c|c|c|c|} \hline 
$\a$   & $\b$  & $f_\a$  & $\c^\a(3,1,1)$  & $\c^\a(2,2,1)$  & $ X^\b $  
  & $F_{\a\b}$ & $ Y^\a$ 
\\ \hline \hline
(5)     &  (4)     &  1  &  1  &  1  &  2  &  3  & 3  \\ \hline
(4,1)   &  (4)     &  4  &  1  &  0  &  2  &  3  & 3  \\ \hline
(4,1)   &  (3,1)   &  4  &  1  &  0  &  2  &  3  & 3  \\ \hline
(3,2)   &  (3,1)   &  5  & -1  &  1  &  2  & -9  & 3  \\ \hline
(3,2)   &  (2,2)   &  5  & -1  &  1  &-10  &  3  & 3  \\ \hline
(3,1,1) &  (3,1)   &  6  &  0  & -2  &  2  &  3  & 3  \\ \hline
(3,1,1) &  (2,1,1) &  6  &  0  & -2  &  2  &  3  & 3  \\ \hline
\end{tabular}
\end{center}
\end{table}

\medno
\textit{Proof of Lemma \ref{thm:w2}}.

\noindent
The proof is a more or less straightforward computation requiring
a little care to avoid double counting.
In order to avoid excessively
cluttered formulas, all variables appearing below the summation symbol
are always summed from $1$ to $n-1$. Moreover, if two or more variable
appear in brackets below the summation, they are assumed to be \textit{all distinct}.
For instance
\begin{equation*}
  \sum_{ j < k,\; [i j \ell] }
  \qquad \text{means}\qquad 
  \sum_{ \atopp{i,j,k,\ell=1}{j < k,\; i\ne j,\; i\ne \ell,\; j \ne \ell} }^{n-1} \,.
\end{equation*}
From \eqref{eq:hw}, letting $A_{ij} := (in) +(jn) -(ij)$, we get
\begin{equation*}
  \hw = 
  \sum_i x_i^2 \, (in) + \ov 2\sum_{[i j]}
  x_i x_j \, A_{ij} \,.
\end{equation*}
Given two elements $a,b$ of the group algebra, we denote with $\{a,b\} = ab + ba$ the
anticommutator of $a$ and $b$. Thus
\begin{equation}
  \label{eq:uvz}
  \hw^2 = \frac{u}{2} + \frac{v}{2} + \frac{z}{8} \,,
\end{equation}
where
\begin{equation}
\begin{split}
  u &:= \sum_{i,j} x_i^2 x_j^2 \, \{(in),(jn)\}
  \hspace{30pt}
  v := \sum_{i,\, [j k]} x_i^2 x_j x_k \, \{ (in), A_{jk} \} 
  \\
  z &:= \sum_{[i j]\, [k \ell]} x_i x_j x_k x_\ell \, \{A_{ij}, A_{k \ell}\} \,.
\end{split}
\end{equation}
If $i,j,k,\ell$ are \textit{distinct} integers in $\{1, \ldots, n-1\}$, we have
the following identities
\begin{align*}
  &\{ (ij), (k\ell)\} = 2 (ij) (k\ell)
  \\
  &\{ (in), (jn)\} = (jnj)+(jni) 
  \\
  & \{(in), A_{jk}\} = (inj) +(jni)+(ink)+(kni)-2(in)(jk)
  \\
  & \{(in), A_{ik}\} = 2 \cdot \uuno_{S_n} 
  \\
  & \{A_{ij}, A_{ij}\} = 2 A_{ij}^2 = 2\bigl( 3\cdot \uuno_{S_n} - (inj) - (jni) \bigr) 
  \\
  & \{A_{ij}, A_{ik}\} =  2\cdot \uuno_{S_n} + (jnk)+(knj)+(ikj)+(jki) \\
  &\hspace{50pt} - 2(ik)(jn) -2 (ij)(kn)
  \,.
\end{align*}  
Let us start with $u$ term. 
\begin{equation}
\begin{split}
  \label{eq:u2}
  u &= 
  2 \sum_{i} x_i^4 \, \uuno_{S_n}
  + 2 \sum_{ [i j] } x_i^2 x_j^2 \, (inj) \,.
\end{split}
\end{equation}  
For $v$ we obtain
\begin{equation}
\begin{split}
  \label{eq:v}
  v &= \sum_{[i j k]} x_i^2 x_j x_k \, \{ (in), A_{jk} \} 
  + 2 \sum_{[j k]} x_j^3 x_k \, \{ (jn), A_{jk} \}  
  \\
  &=
  \sum_{[i j k]} x_i^2 x_j x_k \, \bigl[ (inj) +(jni)+(ink)+(kni)-2(in)(jk) \bigr]
  \\
  &+ 
  4 \sum_{[j k]} x_j^3 x_k \, \uuno_{S_n} \,.
\end{split}
\end{equation}
We split $z$ as a sum $z=z_1+z_2+z_3$, according to the cardinality of 
$\{i,j\}\cap\{k,\ell\}$, as follows
\begin{align*}
  z &:= 2 \sum_{[i j]} x_i^2 x_j^2  \, \{A_{ij}, A_{ij}\} 
  + 4 \sum_{[i j k]} x_i^2 x_j x_k  \, \{A_{ij}, A_{ik}\} \\
  &+ \sum_{[i j k \ell]} x_i x_j x_k x_\ell  \, \{A_{ij}, A_{kl}\} 
  =: z_1 + z_2 + z_3 
  \,.
\end{align*}
We have
\begin{align*}
  z_1 &:= 
  \sum_{[i j]} x_i^2 x_j^2  \, \bigl( 12\cdot \uuno_{S_n} - 8 (inj) \bigr) 
  \,.
\end{align*}
Hence
\begin{align}
  \label{eq:uz1}
  \frac{u}{2} + \frac{z_1}{8} =   \sum_{i} x_i^4 \, \uuno_{S_n} +
  3 \sum_{i<j} x_i^2 x_j^2  \,  \uuno_{S_n} \,.
\end{align}
Next, using the anticommutation relations we found above, we find
\begin{equation}
\begin{split}
  \label{eq:vz2}
  \frac{v}{2} + \frac{z_2}{8} &:= 
  \ov2 \sum_{[i j k]} x_i^2 x_j x_k \bigl[ 2\cdot \uuno_{S_n} + J^{(3,1)}_{\{i,j,k,n\}}
  - 2 J^{(2,2)}_{\{i,j,k,n\}} \bigr] 
  + 
  4 \sum_{[j k]} x_j^3 x_k \, \uuno_{S_n} 
  \\
  &=
  \ov2 \sum_{[i j k]} x_i^2 x_j x_k \bigl[ 2\cdot \uuno_{S_n} + X_{\{i,j,k,n\}} \bigr] 
  + 
  4 \sum_{[j k]} x_j^3 x_k \, \uuno_{S_n} 
  \,.
\end{split}
\end{equation}  
Thanks to the symmetry of the four indices, the term $z_3$ can be written as
\begin{equation}
\begin{split}
  \label{eq:z3}
  z_3 &= 
  2 \sum_{[i j k \ell]} x_i x_j x_k x_\ell  \,  A_{ij} A_{kl} \\
  &=
  \frac{2}{3} \sum_{[i j k \ell]} x_i x_j x_k x_\ell  \,  
  ( A_{ij} A_{kl} + A_{ik} A_{j\ell} + A_{i\ell} A_{jk} ) \\
  &=
  16 \sum_{i < j < k < \ell} x_i x_j x_k x_\ell  \,  
  ( A_{ij} A_{kl} + A_{ik} A_{j\ell} + A_{i\ell} A_{jk} ) \,.
\end{split}
\end{equation}
The term $A_{ij} A_{kl} + A_{ik} A_{j\ell} + A_{i\ell} A_{jk}$ contains: (a) all
3-cycles $(inj)$, with $1 \le i,j \le n-1$, with coefficient $+1$; 
(b) all double transpositions $(in)(jk)$ with $1 \le i,j,k \le n-1$, with coefficient $-1$;
(c) all double transpositions $(ij)(k\ell)$ with $1 \le i,j,k,\ell \le n-1$, with coefficient $+1$.
Hence we have
\begin{equation}
\begin{split}
  \label{eq:A}
  A_{ij} A_{kl} + A_{ik} A_{j\ell} + A_{i\ell} A_{jk} &= 
  J^{(3,1,1)}_{\{i,j,k,\ell,n\}} - J^{(3,1)}_{\{i,j,k,\ell\}}
  - J^{(2,2,1)}_{\{i,j,k,\ell,n\}} + 2 J^{(2,2)}_{\{i,j,k,\ell\}} \\
  &= Y^{n}_{\{i,j,k,\ell\}} \,.
\end{split}
\end{equation}
The proof of the lemma now follows 
from \eqref{eq:uvz}, \eqref{eq:v}, \eqref{eq:uz1}, \eqref{eq:vz2}, \eqref{eq:z3} and \eqref{eq:A}.
\qed

\medno

%--------------------------------------------------------------------
\section{The Kazhdan constant and (im)possible generalizations} 
\label{sec:kazh}
%--------------------------------------------------------------------

\noindent
Let $G$ be a finite group, let $G^*$ be the set of all 
(not necessarily finite-dimensional) unitary representations
of $G$ and let $G^*_0$ be the set of all $\dR \in G^*$ with no
invariant nonzero vector, that is
\begin{align*}
  G_0^* := \{ \dR \in G^* : \dI \not\subset \dR \} \,.
\end{align*}
If $(\dR,V) \in G_0^*$ and $Q$ is a generating set of $G$,
the \textit{Kazhdan constant} associated to the pair $(Q, \dR)$
is defined as
\begin{align}
    \label{eq:kazmin}
    \k_{G}(Q,\dR) := \inf_{v\in V : \|v\|=1} \max_{q\in Q} \| \dR(q) v -v \| \,.
\end{align}
By minimizing over $\dR \in G_0^*$ we obtain the \textit{Kazhdan constant of $Q$}
\begin{align*}
  \k_{G}(Q) := \inf_{\dR\in G_0^*} \k_{G}(Q,\dR) \,.
\end{align*}
We can always assume that $Q$ is a symmetric generating set, that is $Q = Q^{-1}$,
since adding the inverse of an element to $Q$ does not change the value
of the Kazhdan constant.
The Kazhdan constant can be defined in the much more general framework
of a topological group $G$. We refer the reader to \cite{BeHaVa} for a review
on on the Kazhdan constant and Kazhdan's property (T) in this more general setting.
We have observed in section \ref{sec:notation} that, being $\psi_G(w)$
a spectral quantity, in its definition we can either minimize over
all representations or just the irreducible ones getting the same
result. This is false for the Kazhdan constant.
If we define
\begin{align*}
  \k^I_G(Q) := \inf_{\dR\in G_0^* \cap \Irr(G) } \k_{G}(Q,\dR) \,,
\end{align*}
then we have \cite{BacHar}
\begin{equation*}
  |Q|^{-1/2} \k^I_G(Q) \le \k_{G}(Q) \le \k^I_G(Q) \,,
\end{equation*}
and explicit examples where the second inequality is strict
are presented in \cite{BacHar}.

In order to illustrate the connection between the Kazhdan constant 
and the spectral gap $\psi$, we start by introducing a slight
abuse of notation, by identifying the  generating set $Q$ with
the group algebra element $\hQ := \sum_{q\in Q} q$, and we observe that
if $Q$ is symmetric, then $\hQ \in \bC G^{(+)}$, 
so (\ref{it:3}) of Proposition \ref{thm:lapl} applies to $\hQ$.
Thus if $\dI \not\subset \dR$, then the associated representation
Laplacian $\D_G(Q,\dR)$ has no trivial eigenvalues.
By consequence, in accord with \eqref{eq:sgdef}, we have
\begin{equation}
  \label{eq:psiYS}
  \psi_G(Q,\dR) 
  := \psi_G(\hQ,\dR ) 
  = \min \spec \D_G(Q,\dR)
  = \min \spec \sum_{q\in Q} (\id - \dR(q) ) \,.
\end{equation}
From \eqref{eq:pos} it follow that if $(\dR,V) \in G_0^*$,
then
\begin{equation}
  \label{eq:lapl=}
  \< \D_G(Q,\dR) v, v \> = \ov2 \sum_{q\in Q} \|\dR(q) v -v\|^2 
  \qquad
  \text{for each $v\in V$.}
\end{equation} 
Since $\D_G(Q,\dR)$ is self-adjoint, we get
\begin{equation*}
  \min \spec \D_G(Q,\dR)  =\inf_{v\in V : \|v\|=1}  \< \D_G(Q,\dR) v, v \> \,.
\end{equation*}
In this way we obtain a formula for the spectral gap which resembles
the definition of the Kazhdan constant
\begin{equation}
    \label{eq:gapmin}
  \psi_G(Q,\dR) = \ov2 \inf_{v\in V: \|v\|=1} \sum_{q\in Q} \|\dR(q) v - v\|^2 \,.
\end{equation}
Comparing with \eqref{eq:kazmin}, we get, for any $\dR \in G_0^*$,
\begin{align}
    \label{eq:kp1}
    \k^2_{G}(Q,\dR) &\le 2 \psi_G(Q,\dR) \le |Q|\, \k^2_{G}(Q,\dR) \\
    \label{eq:kp2}
    \k^2_{G}(Q)     &\le 2 \psi_G(Q) \le |Q|\, \k^2_{G}(Q) \,.
\end{align}
In \cite{Kas} Kassabov, using an approach
completely different from \cite{CaLiRi} 
proved that the ``minimality'' property  \eqref{eq:D} of the defining
representation for computing the spectral gap \textit{and the Kazhdan constant}
also holds for finite Coxeter systems $(G,Q)$.

A Coxeter group  $G$ generated by $Q = \{s_1, \ldots , s_n\}$ 
is defined  by integral numbers $m_{\a \b} = m_{\b\a}$ such that
$m_{\a\a}=1$ and $m_{\a\b}\ge 2$ if $\a\ne\b$. 
$G$ has presentation (see \cite{Hum} for background)
\begin{equation*}
  G \cong \bigl\< s_\a \tc \ (s_\a s_\b)^{m_{\a \b}} = 1, 
  \forall\ \a, \b  \bigr\> \,.
\end{equation*}
It is known that $G$ has a \textit{defining representation} 
on a $n$-dimensional vector space $V$ such that
each generator $s_\a$ acts as a reflection with respect to a hyperplane $V_\a$. 
Moreover if $G$ is finite the angle between the 
hyperplanes $V_\a$ and $V_\b$ is equal to $\pi/m_{\a\b}$.

\smallno
In \cite[Thm. 1.3]{Kas} the following result has been proved:

\begin{theorem}\label{thm:kas}
Let $G$ be a finite Coxeter group with Coxeter generating set $Q$.
Then 
\begin{align}
  \label{eq:kas}
  &\psi_G(Q) = \psi_G(Q,\dD')
  &\k_G(Q) = \k_G(Q, \dD')
\end{align}
where $\dD'$ is the defining representation of $G$.
\end{theorem}

\begin{HNA}\label{thm:amb}
When dealing with the symmetric group $S_n$ 
one usually call defining representation the \textit{reducible} $n$-dimensional
representation $\dD$ associated with the natural action of
$S_n$ on the set $\{1,2,\ldots,n\}$, which decomposes as in \eqref{eq:DI}.
On the other side,
$S_n$, when generated
by $n-1$ adjacent transpositions
$Q = \{(1,2),\, (2,3),\, \dots,\, (n-1,n) \}$ corresponds
to Coxeter group $A_{n-1}$ and, in this context, 
its defining representation
is the $(n-1)$-dimensional \textit{irreducible} representation 
$\dD' \cong [n-1,1]$. This is really harmless, since 
$\psi_G(Q,\dD) = \psi_G(Q,\dD')$, so the left half of \eqref{eq:kas}
coincides with \eqref{eq:D}.
\end{HNA}

\medno
Kassabov's theorem is more general in some respects and less
general in others, than Theorem \ref{thm:ald}.
It is less general because it applies (essentially) to
only one element of the group algebra, namely $w = \hQ = \sum_{q\in Q} q$,
where $Q$ is a Coxeter generator for $G$.
Thus, for instance, in the case of the symmetric group,
it holds for $\psi_{S_n}(w)$ when
$w$ is  a sum of distinct adjacent transpositions,
hence (modulo permutations of the indices) 
\begin{equation*}
  w = (1,2) + (2,3) + \cdots + (n-1,n) \,.
\end{equation*}
It is more general because it holds for any finite Coxeter group
and it applies to the Kazhdan constant,
These considerations prompt the surfacing of at least two questions:
\begin{enumerate}[(1)]
\item 
Can Kassabov's approach be used for dealing with more
general elements $w$ of the group algebra or more general generating sets?
\item
Does the analogous of Theorem \ref{thm:ald} hold for the Kazhdan
constant?
\end{enumerate}
The answers are (1) not likely and (2) no.

\medno
The first answer can be explained as follows: 
a key ingredient in the proof of Theorem \ref{thm:kas} is that the 
symmetric matrix
$A$ with elements
\begin{equation}
  \label{eq:Aab}  
  A_{\a\b} = - \cos \frac{\pi}{m_{\a\b}}
\end{equation}
is (strictly) positive definite as a quadratic form.
The problem is that positive definite matrices of the form \eqref{eq:Aab}
have been classified (see, for instance, \cite{Hum}) and 
correspond exactly to finite Coxeter systems $(G,Q)$, so there is little
hope to go beyond Coxeter systems with this approach.
To consider a concrete example, let us take $S_4$ with \textit{non-Coxeter} generating set
\begin{equation*}
  Q = \{ s_\a=(1,2),\, s_\b=(1,3),\, s_\g=(1,4) \} \,.
\end{equation*}
Then it easy to check that 
$(s_\mu s_\nu)^{m_{\mu\nu}}=1$ where $m_{\a\b} = m_{\a\g} = m_{\b\g} = 3$.
Thus we have
\begin{equation*}
  A = 
  \begin{pmatrix}
      1 & -1/2 & -1/2 \\ -1/2 & 1 & -1/2 \\ -1/2 & -1/2 & 1 
  \end{pmatrix}
  \qquad
  \det A = 0 \,,
\end{equation*}
hence $A$ is positive semidefinite, but not positive definite
and the approach of Theorem \ref{thm:kas} will not work in this case,
even if the left half of \eqref{eq:kas} holds, 
thanks to Theorem \ref{thm:ald}.

\medno
We now present a family of counterexamples which justify
our negative answer to question (2).

\begin{proposition}\label{thm:ce}
Let $G=S_n$ with $n\ge 4$ and let $T_n$ be the set of all transpositions
of $S_n$. Then, if $\dD' \cong [n-1,1]$ is the defining representation,
we have
\begin{align*}
  \k_{S_n}( T_n, \dD') < \k_{S_n}( T_n) = \frac{2}{\sqrt{n-1}} \,.
\end{align*}
\end{proposition}

\smallno
We start with a couple of preliminary results.

\begin{lemma}\label{thm:compl}
For each integer $n\ge 2$ we have $\psi_{S_n}(T_n) = n$.
\end{lemma}

\Pro\ 
This is pretty much well known. It can be easily deduced, for instance,
from \cite[Cor. 4]{DiSh}, taking into account their different 
normalization. Since the proof is a simple combination of facts 
previously discussed in this paper, we include it for the sake of completeness.
From Theorem \ref{thm:ald}, \eqref{eq:IT} and Lemma \ref{thm:schur}
applied to the conjugacy class of transpositions with cycle partition 
$\a = (2,1^{n-2})$, we get
\begin{equation}
\begin{split}
    \label{eq:gap1}
  \psi_{S_n}(T_n) &= \psi_{S_n}(T_n, \dD) = \psi_{S_n}(T_n, [n-1,1]) 
  \\
  &= |T_n| - \l^*\bigl[ \dT^{(n-1,1)}( J^{(2,1^{n-2})} ) \bigr]
  \\
  &=
  |T_n| - \frac{ \c^{(n-1,1)}( (2,1^{n-2}) ) \, |T_n| } { f_{(n-1,1)} } \,.
\end{split}
\end{equation}  
In the case of transpositions,
Frobenius formulas for the irreducible characters
take the simple form \cite{Ing}
\begin{equation}
  \label{eq:q=}
  \frac{ \c^{\b}( (2,1^{n-2}) ) \, |T_n| } { f_{\b} } 
  = \ov2 \sum_{i=1}^r \b_i \, [ \b_i - (2i-1)] \,,
\end{equation}
where $r$ is the \textit{length} (number of terms) of the partition $\b$.
The Lemma follows form \eqref{eq:gap1} end \eqref{eq:q=} with $\b=(n-1,1)$.
\qed

\medno
The second result shows that if 
the second
inequality in \eqref{eq:kp1} is saturated, then we are in a 
very special situation.

\begin{lemma}\label{thm:spe}
Let $G$ be a finite group. If $Q$ is a generating set for $G$ 
and $(\dR,V) \in G_0^*$ is finite-dimensional and such that 
\begin{equation}
  \label{eq:satu}
  |Q| \, \k^2_{G}(Q,\dR) = 2 \psi_G(Q, \dR)   \,,
\end{equation}
then there exists a unit vector $u\in V$ such that
\begin{equation}
  \label{eq:alls}
  \k_{G}(Q,\dR) = \| \dR(q) u - u \| 
  %% = \frac{2}{|Q|} \psi_G(Q, \dR)
  \qquad
  \forall q\in Q\,.
\end{equation}
\end{lemma}

\Pro\ 
Since $\dR$ is finite-dimensional, in \eqref{eq:kazmin}
we are minimizing a continuous function on a compact set, hence
there exists a unit vector $u \in V$ such that 
\begin{align}
  \label{eq:infmin}
  \max_{q\in Q} \|\dR(q) u - u \| =
  \k
_G(Q,\dR) 
 \,.
\end{align}
If \eqref{eq:satu} holds, we have
\begin{equation}
\begin{split}
    2 \psi_G(Q,\dR) &\le \sum_{q\in Q} \|\dR(q) u - u \|^2
    \le |Q| \, \max_{q\in Q} \|\dR(q) u - u \|^2
    \\
    &= |Q| \, \k^2_G(Q,\dR) = 2 \psi_G(Q, \dR) \,.
\end{split}
\end{equation}  
Hence the second (as well as the first) inequality above must be an
equality, which implies the Lemma.
\qed

\medno
\textit{Proof of Proposition \ref{thm:ce}}. 
Since $T_n$ is a conjugacy class, the group of all inner automorphisms
of $S_n$ leaves $T_n$ invariant and acts transitively on $T_n$. 
Hence we can apply Proposition 1 of \cite{PakZuk}
which, in our notation, states that 
the second inequality in \eqref{eq:kp2} is saturated. 
Thus, using Lemma \ref{thm:compl}, 
we have
\begin{equation}
    \label{eq:satu1}
  \k^2_{S_n}(T_n) = \frac{2}{|T_n|} \psi_{S_n}(T_n) 
  = \frac{4}{n-1} \,.
\end{equation}
We intend to show that the assumption 
\begin{equation}
  \label{eq:abs}
  \k_{S_n}(T_n) = \k_{S_n}(T_n,\dD')
\end{equation}
leads to a contradiction. In fact by Theorem \ref{thm:ald} we have
$\psi_{S_n}(T_n) = \psi_{S_n}(T_n, \dD')$, so, if \eqref{eq:abs}
holds, thanks to \eqref{eq:satu1} we can apply Lemma \ref{thm:spe}
to the pair $(T_n, \dD')$ and we could conclude that, for some
unit vector $u$ in the representation space of $\dD'$, we have
\begin{equation}
  \label{eq:alleq}
  \| \dD'(q) u - u \|^2 = \k^2_G(T_n, \dD')
  = \frac{4}{n-1}
  \qquad
  \forall q\in T_n\,.
\end{equation}
But this eventuality can be easily ruled out. 
Let $V$ be an $n$-dimensional complex vector space with
orthonormal basis $(e_k)_{k=1}^n$,
and let $\dD$ be the defining $n$-dimensional 
representation of $S_n$ acting on $V$ as
\begin{equation*}
  \dD( \pi )\( \sum_{k=1}^n u_k e_k \) = 
  \sum_{k=1}^n u_k e_{\pi(k)}
  \qquad
  u_k \in \bC,\ \pi \in S_n \,. 
\end{equation*}
The representation $\dD'$ can be realized by restricting $\dD$
to the invariant subspace 
\begin{equation*}
V' := \Bigl\{ v \in V : \Bigl\< v, \sum_{k=1}^n e_k \Bigr\> =0 \Bigr\} \,.
\end{equation*}
Let then $u = \sum_{k=1}^n u_k e_k$ be a unit vector in $V'$ such that
\eqref{eq:alleq} holds. Then, by explicit computation, for each transposition
$(i,j)$, we have
\begin{equation*}
  \| \dD'((i,j)) u - u \|^2 = 2 |u_i-u_j|^2
  = \frac{4}{n-1}
  \qquad
  \forall i\ne j \,.
\end{equation*}
This implies, in particular, that there are $n$ points
$u_1,\ldots, u_n$ in the complex plane which are equidistant
from each other, which is impossible if $n\ge 4$.
\qed

\begin{remark}\label{thm:rem1}
The condition $n\ge 4$ cannot be omitted. 
For $n=3$ we have in fact $\k_{S_n}(T_n) = \k_{S_n}(T_n, \dD')$.
Proceed as in the previous proof and take 
$u=(1, \nep{i2\pi/3}, \nep{-i2\pi/3})/\sqrt{3}$.
\end{remark}

\begin{remark}\label{thm:rem2}
A minimizing representation $\dR$, for which we have the equality 
$\k_{S_n}(T_n,\dR) = \k_{S_n}(T_n)$
is actually a direct sum of a certain number of copies of $(\dD',V)$.
Proceeding more or less as in the proof of Proposition 1 of \cite{PakZuk}
(see also the proof of Theorem 1.2 of \cite{Neu}),
let 
\begin{equation*}
    \dR := |T_n| \cdot \dD' = \dD' \oplus \cdots \oplus \dD' 
    \qquad\text{($|T_n|$ terms)}.
\end{equation*}
Let $u \in V$ be a unit minimizing vector in \eqref{eq:gapmin} for $\psi_{S_n}(T_n,\dD')$,
so that
\begin{equation*}
  \psi_{S_n}(T_n, \dD') = \frac{1}{2} \sum_{t\in T_n} \|\dD'(t) u - u \|^2 \,.
\end{equation*}
Let
\begin{equation*}
  \hhU = (\dD'(t) u)_{t\in T_n} \in \hV := V \oplus \cdots \oplus V \,.
\end{equation*}
Then, since $T_n$ is a conjugacy class and $\dD'$ is unitary, we get
\begin{align*}
  \| \dR(q) \hhU - \hhU\|^2 
  &= \sum_{t\in T_n} \| \dD'(q t) u - \dD'(t)u \|^2
  = \sum_{t\in T_n} \| \dD'(t^{-1} q t) u - u \|^2
  \\
  &= \sum_{t\in T_n} \| \dD'(t) u - u \|^2 \,.
\end{align*}
Hence the LHS is independent of $q$, which implies
\begin{align*}
  \k^2_{S_n}( T_n, \dR) 
  &\le \frac{  \| \dR(q) \hhU - \hhU\|^2 }{ \|\hhU\|^2 }
  = \ov{|T_n|} \sum_{t\in T_n} \| \dD'(t) u - u \|^2
  \\
  &= \frac{2}{|T_n|} \psi_{S_n}(T_n, \dD') = \frac{4}{n-1} \,.
\end{align*}
On the other side, by \eqref{eq:satu1}, we know that 
$\k^2_{S_n}(T_n) = 4/(n-1)$, hence we have
$\k_{S_n}(T_n,\dR) = \k_{S_n}(T_n)$.
\end{remark}

\medno
We conclude with a few words about the possibility of generalizing
Theorem \ref{thm:ald} to elements $w$ of the group algebra
which are no longer combinations of just transpositions.
In \cite{Ces3} we have proved \eqref{eq:D} when $w$ is the sum of all
initial reversals
\begin{equation*}
  w = r_1 + r_2 + \cdots + r_n \,,
\end{equation*}
where $r_k \in S_n$ is the permutation which reverses the order 
of the first $k$ positive integers.
General results however do not appear easy to obtain.
Even in very special (and simple) case of $w = J^\a$
which is the sum of all elements of the conjugacy class $C^\a$
(remember \eqref{eq:Ja}), identity \eqref{eq:D} in general fails.
If $\a$ is a class of even permutations, the conjugacy class $C^\a$
does not generate $S_n$ and we have, trivially
\begin{equation*}
  \psi_{S_n}(J^\a) = \psi_{S_n}(J^\a, [1^n]) = 0 \,.
\end{equation*}
A counterexample for an odd class is given by 
$\a = (4,1)$, 
the set of all $4$-cycles
in $S_5$. It follows from Lemma \ref{thm:schur} that
\begin{equation*}
  \psi_{S_5}( J^{(4,1)} ) = \psi_{S_5}( J^{(4,1)}, [2,2,1] ) = 24 
  < 30 = \psi_{S_5}( J^{(4,1)}, \dD ) \,.
\end{equation*}
There is one conjecture though, by P. Caputo, which is defying our
attempts to either prove it, or disprove it with 
numerical experiments.
For each $A \sset \{1,\ldots, n\}$, 
let $S_{n,A}$ be 
the set of all $\pi \in S_n$ such that $\pi(i) = i$ for each $i \in A^c$, and let
\begin{equation*}
  J_{n,A} := \sum_{\pi \in S_A} \pi \,.
\end{equation*}
In other words $J_{n,A}$ is the sum of all possible shuffles
inside $A$ which leave each element outside invariant.
The conjecture is as follows:

\begin{capu}\label{thm:capu}
(P. Caputo, private communication).
Let $n$ be a positive integer, $n\ge 3$. 
For each $A \sset \{1,\ldots, n\}$, let $\a_A \ge 0$ and let
\begin{equation*}
  w = \sum_{A \sset \{1, \ldots, n\}}  \a_A \, J_{n,A}.
\end{equation*} 
If $\supp w$ generates $S_n$, then $\psi_{S_n}(w) = \psi_{S_n}(w, \dD)$.
\end{capu}

\medno
This may be considered a rather natural generalization of Theorem \ref{thm:ald},
since, by imposing the extra condition $|A|=2$, one recovers
Aldous' conjecture.

%--------------------------------------------------------------------
%--------------------------------------------------------------------

\noindent
\bibliographystyle{amsalpha-op}
%\bibliographystyle{amsplain-op}
%
% do: convert-bibtex.pl < ref.bib >ref-no-some-fields.bib
%\bibliography{ref-no-some-fields}

\def\polhk#1{\setbox0=\hbox{#1}{\ooalign{\hidewidth
  \lower1.5ex\hbox{`}\hidewidth\crcr\unhbox0}}}
\providecommand{\bysame}{\leavevmode\hbox to3em{\hrulefill}\thinspace}
\providecommand{\MR}{\relax\ifhmode\unskip\space\fi MR }
\providecommand{\MRhref}[2]{%
  \href{http://www.ams.org/mathscinet-getitem?mr=#1}{#2}
}
\providecommand{\href}[2]{#2}

%------------------------------------------------------------------------
\end{document}